# An Efficient Algorithm to Calculate the Center of the Biggest Inscribed Circle in an Irregular Polygon

OSCAR MARTINEZ

*A*BSTRACT   *In this paper, an efficient algorithm to find the center of the biggest circle inscribed in a given polygon is described. This work was inspired by the publication of Daniel Garcia-Castellanos & Umberto Lombardo and their algorithm [1] used to find a landmass' poles of inaccessibility. Two more efficient algorithms were found, one of them only applicable when the problem can be described as a linear problem, like in the case of a convex polygon.*

**Introduction**

Given a polygon of *N* sides, how can we find the biggest inscribed circle? In other words, what is the furthest point from any vertex or edge? If the polygon has 3 sides, which would make it a triangle, this point is equivalent to the incircle and could be found by intersecting the bisection of any two of its angles. Likewise, if the polygon is a square, a pentagon, or any other regular shape, the same method works as illustrated in Figure 1. However, when that method is applied to an irregular polygon with more than 3 sides it might fail as shown in Figure 2.

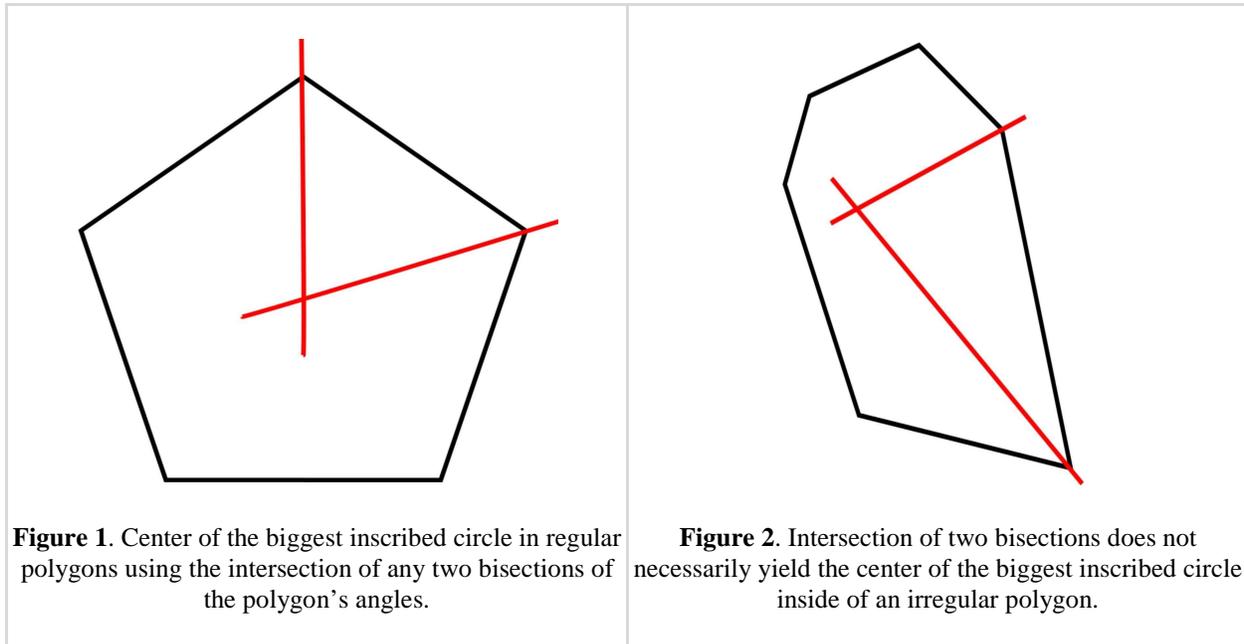

**Figure 1**. Center of the biggest inscribed circle in regular polygons using the intersection of any two bisections of the polygon's angles.

**Figure 2**. Intersection of two bisections does not necessarily yield the center of the biggest inscribed circle inside of an irregular polygon.

When facing the problem with an irregular polygon, the solution can be found by taking the same approach used for the Poles of Inaccessibility (PIA) problem. A PIA is the most difficult point to reach within a landmass, commonly described as the point furthest from any coastline. If we look at a landmass as a polygon, the PIA of such landmass is identical to the center of the biggest circle inscribed in it. Summarizing the known algorithm's methodology, to find the PIA we first describe a region *R* around a initial candidate for the PIA, then divide *R* into *n* by *m* cells and find the node (intersection between lines) furthest from any coastline to make it our next candidate PIA. Subsequently, *R* is centered around the new candidate PIA and shrunk by a factor of *k*, finding a new candidate and repeating the process until the desired precision is achieved. Figure 3 shows a visual demonstration of the methodology; for a detailed description of how the sequential algorithm used to find the PIA works refer to the original paper by Daniel Garcia-Castellanos & Umberto Lombardo.

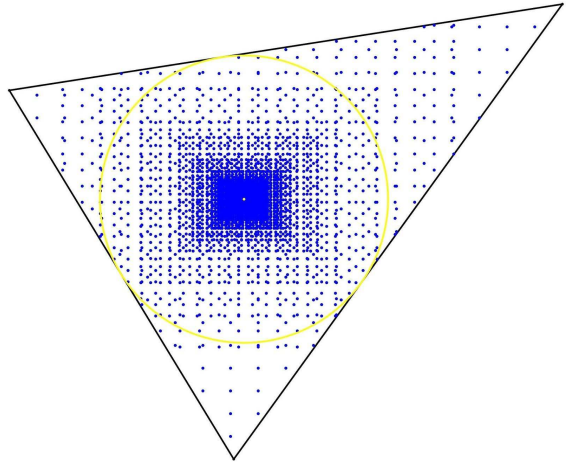

**Figure 3**. Visual demonstration of the sequential algorithm to find the PIA. Each point represents a node that has been compared against the candidate PIA in its iteration and the algorithm is recursively called until the density of nodes corresponds to the desired precision to be achieved.

## Monte Carlo's approach

The first proposed methodology is very similar to the known algorithm used to find the PIA. The key difference is in the choice of nodes, which would be done by the use of pseudo-randomly chosen coordinates within the region *R*, rather than sequentially by intersecting *n* and *m* lines. By choosing the nodes pseudo-randomly, the an iteration of the algorithm may finish at any arbitrary point rather than having to test $n \times m$ nodes in each of the regions. For this implementation, the region is shrunk when a new candidate PIA has not been found within the last *k* attempts. In Figure 4 we can visualize how the algorithm works and in Figure 5 the pseudo-code describing the algorithm is presented.

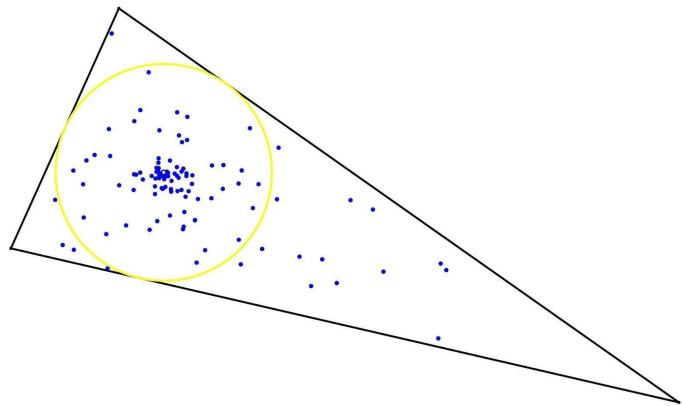

**Figure 4**. Visual demonstration of the randomized algorithm to find the PIA. Each point represents a node that has been compared against the candidate PIA in its iteration and the algorithm is recursively called until the density of nodes corresponds to the desired precision to be achieved.

```
while accuracy > minimum_accuracy do:

    # begin loop through nodes
    while consecutive_misses < k do:

        # select coordinates at random within bounds
        x := random(min_x, max_x)
        y := random(min_y, max_y)
        node := (x, y)

        # loop through edges and find shortest distance
        for edge := first_edge to last_edge do:
            if distance(node, edge) < smallest_distance
                smallest_distance := distance(node, edge)
                PIA := node
            end if
        end do

        # maximize the minimum distance through iterations
        # and keep track of the consecutive number of times
        # that a smallest distance has not been found
        if smallest_distance > maximin_distance
            maximin_distance := smallest_distance
            consecutive_misses := 0
        else
            consecutive_misses := consecutive_misses + 1
        end if
    end do

    # calculate current level of accuracy based on the
    # smallest distance between the upper and lower bound
    accuracy := min(max_x - min_x, max_y - min_y)

    # update the bounds of the region
    min_x := PIA(x) - (max_x - min_x) / (sqrt(2) * 2)
    max_x := PIA(x) + (max_x - min_x) / (sqrt(2) * 2)
    min_y := PIA(y) - (max_y - min_y) / (sqrt(2) * 2)
    max_y := PIA(y) + (max_y - min_y) / (sqrt(2) * 2)

end do

return PIA
```

**Figure 5**. Pseudocode describing the algorithm using a randomized choice of coordinates.

Although the big O of the proposed methodology is identical to the known sequential approach, as long as we take as many random coordinates as nodes in the sequential algorithm, it is expected, yet not guaranteed, to find a solution with the same level of accuracy significantly quicker by using a smaller number of nodes. The caveat of this methodology is that, assuming the pseudo-random number generator of choice is effectively random[1], the algorithm is by definition non-deterministic, that is, it could find a solution to the same problem in a significantly different time when used more than once or it could even not find the solution at all! However, if the parameters are well-chosen, the algorithm is almost guaranteed to yield a solution in less or equal time than the known algorithm.

**Linear programming**

If the problem can be described linearly, it is possible to use the well-known linear programming (LP) techniques to solve it in a fraction of the time compared to the two previously discussed algorithms. In order to describe the problem in a form suitable for linear programming, we could try to maximize the minimum distance of the goal point to all of the edges making it a maximin problem. As seen in Figure 6, the set of equations describing the problem is the following:

---

Constraints:

$$d_{\lambda,\varphi}^{\lambda_i,\varphi_i} = \arccos(\sin(\varphi_i) \cdot \sin(\varphi) + \cos(\varphi_i) \cdot \cos(\varphi) \cdot \cos(\lambda_i - \lambda)) \geq Z$$

Objective:

$$\textit{maximize } Z$$

**Figure 6**. Set of constraints and objective function that, once optimized, would yield the solution to the problem. The lack of linearity of some of the equations prevents the use of LP to solve the problem. The point ( $\lambda_i$ , $\varphi_i$ ) represents the closest point to ( $\lambda$ , $\varphi$ ) for the *i*th edge, so $i = 1, \ldots, N$.

---

Because the original problem was characterized by using the Earth's coordinates as the set of all possible solutions, which are approximately spherical coordinates, the equation used to find the great-circle distance between two points is derived from the spherical trigonometry cosine rule. Note that the variable Z is used as the decision variable based on Wald's maximin model[3], so it will become the only variable that we must optimize for in the objective function.

It might be possible to linearize the spherical distance function for optimization purposes but such linearization, if achievable is neither trivial nor apparent. As a result, the proposed problem in Figure 5 could only be solved by the use of non-linear programming. However, if we describe the problem to be in the euclidean plane instead of using spherical coordinates, or $\mathbb{R}^2$, we would have the following set of constraints with the same objective function:

---

[1] A pseudo-random number generator that could output truly random numbers is technically impossible to achieve with classical computer architecture, but it is assumed that its output is *random enough* for practical purposes.

```
Constraints:
```
$$d_{\lambda,\varphi}^{\lambda_i,\varphi_i} = \sqrt{(x_i - x)^2 + (y_i - y)^2} \geq Z, \qquad i = 1, \ldots, N$$

```
Objective:
    maximize Z
```

**Figure 7**. Set of constraints and objective function for the same problem in the Euclidean plane. The point ( $x_i$ , $y_i$ ) represents the closest point to ( $x$ , $y$ ) for the *i*th edge.

Although the problem expressed in Figure 7 is not linear either, we can derive [2] the distance between a point and a line into a more suitable form for linear programming. We would then have the following set of constraints with, once again, the same objective function:

```
Constraints:
```
$$d_{x,y}^{y=m_j x + n_j} = \frac{1}{\sqrt{m_j^2 + 1}} \cdot (y - m_j x + n_j) \geq Z, \qquad j = 1, \ldots, N$$

$$d_{x,y}^{y=m_k x + n_k} = \frac{-1}{\sqrt{m_k^2 + 1}} \cdot (y - m_k x + n_k) \geq Z, \qquad k = 1, \ldots, N$$

```
Objective:
    maximize Z
```

**Figure 8**. Linearization of the problem describing the center of the biggest circle inscribed in a convex polygon in $\mathbb{R}^2$. Note that the set of lines in $y = m_k x + n_k$ are the top edges of the polygon and the set of lines in $y = m_j x + n_j$ are the bottom edges of the polygon, thus $J + K = N$.

Using any linear programming solver, we can solve the problem described by the constraints and objective function in Figure 8 to find the center of the biggest inscribed circle. For this paper, a linear programming solver based on the simplex method was implemented from scratch in Ansi C. The code used to generate the polygons and to formulate the problem is also in C and everything is available at https://github.com/omtinez/CenterPolygon under a permissive BSD license.

Unfortunately, by using the distance from the point to a line as our constraints, we limit the set of polygons suitable to use this algorithm to convex polygons. In order to test the algorithms, three points were chosen pseudo-randomly to form a triangle and the speed and accuracy of the described algorithms were measured.

## Results

Using pseudo-random coordinates, 200 triangles were generated. Out of the three proposed algorithms, the only one that found the exact solution[2] in all of the instances was the algorithm using the linear programming formulation. The charts below summarize the results:

| Results | Sequential N = M = 12 | Sequential N = M = 15 | Sequential N = M = 20 | Randomized K = 15 | Randomized K = 30 | Randomized K = 50 | LP |
|---|---|---|---|---|---|---|---|
| **Exact solution (%)** | 76.33 | 77.51 | 78.1 | 56.27 | 61.06 | 69.40 | 100 |
| **Error <= 0.01% (%)** | 14.79 | 15.38 | 17.75 | 27.75 | 29.17 | 25.09 | 0 |
| **Error <= 0.1% (%)** | 1.18 | 2.96 | 0 | 2.54 | 2.84 | 1.72 | 0 |
| **Error > 1% (%)** | 7.69 | 4.14 | 4.14 | 13.43 | 6.92 | 3.78 | 0 |
| **Runtime (seconds)** | 2.065 | 3.2003 | 5.3291 | 0.5293 | 0.6995 | 1.644 | 0.012 |

**Figure 9**. Table summarizing the results after testing the algorithms on 200 pseudo-randomly generated triangles and taking the average of 10 independent tries.

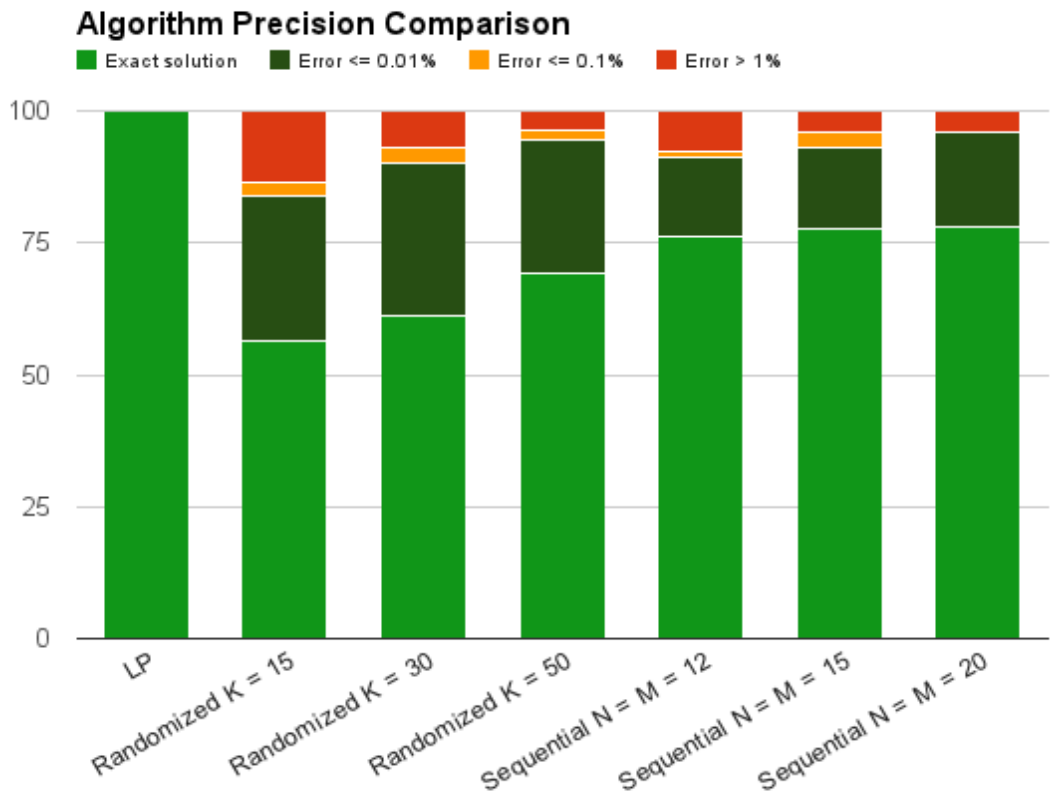

**Figure 10**. Comparison of the average precision yielded by each algorithm.

---

[2] A solution is considered exact if the difference between the solution found by non-numeric methods and the solution yielded by the algorithms is less than the floating point precision of a standard 64-bit Intel CPU, which makes it virtually and indistinguishable for practical purposes.

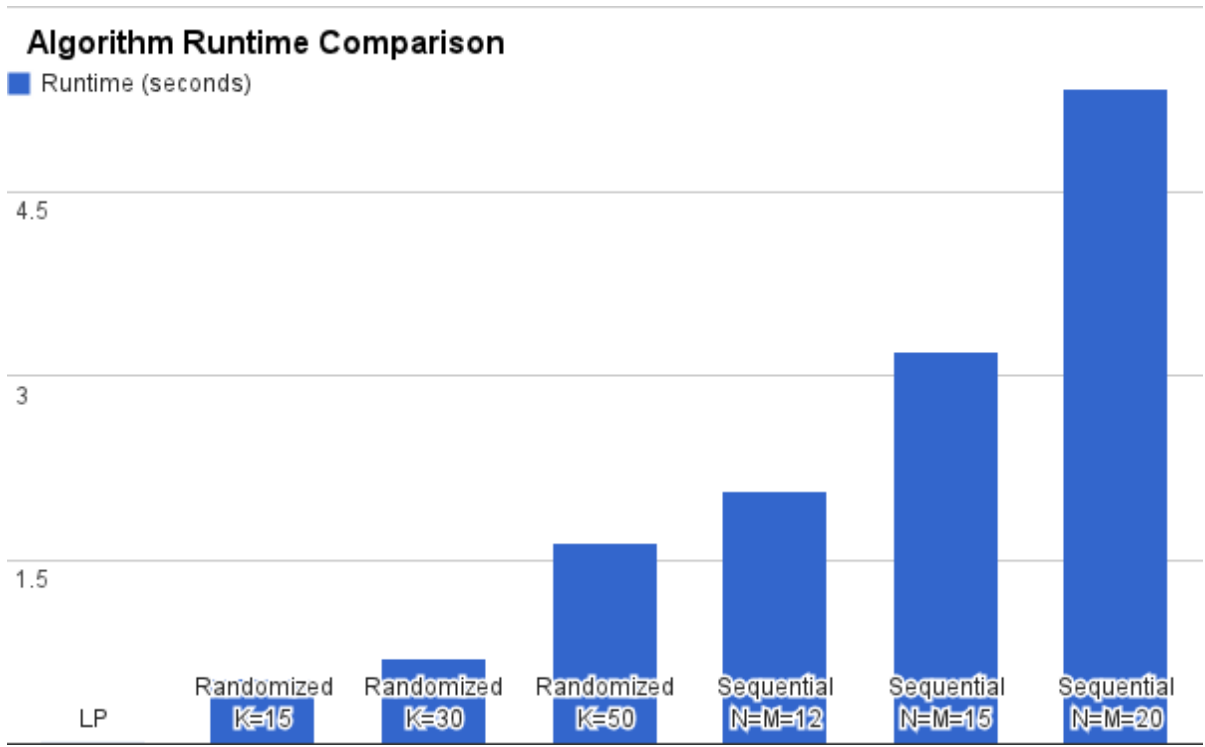

**Figure 11**. Comparison of the average runtime of each algorithm.

From these results, the first and most obvious detail worth noticing is that the sequential and randomized algorithms are very sensitive to the parameters used. Once the number of nodes tested in each iteration is high enough, increasing them further will lead to an increase in runtime at the same rate (almost exponentially) but without yielding significantly better results. Another significant discovery is that both the sequential and the randomized algorithms yield results with an error large enough to make them unsuitable if very high precision is required; a closer look at those instances reveals that the source of those errors were polygons with very small areas like narrow triangles, and the randomized algorithm dealt better with those because in the implementation used, only nodes that lay inside of the polygon count towards the number of nodes used in the iteration. Lastly, the linear programming approach takes virtually no time since most of its runtime is spent in the overhead dedicated to setting up the problem. Because of the way the different algorithms work, it is be expected that using eligible polygons with a larger number of sides would only make a difference in the performance of the sequential and the randomized algorithms, whereas the linear programming algorithm should remain with a near-equal runtime.

**Conclusion**

The clear winner is the linear programming technique, vastly over performing the other two algorithms, but unfortunately it will only work on a limited type of polygons. If the polygons to be dealt with are known to be of such type then it should be the algorithm of choice, but if we are looking for a more general purpose solution the randomized method is the best way to go. Leaving the performance in terms of time and precision aside, the randomized method has the advantage of being relatively friendlier to write an implementation for and has no inherent nested loops, although they can be avoided with the clever use of the *mod* operator in the implementation of the sequential algorithm. If the lack of certainty when looking for the solution is a concern due to the use of random number generators, the randomized algorithm could be combined with a genetic algorithm to decrease the chances of having "bad luck" with the random number generator.

**References**

[1] Garcia-Castellanos, Daniel, and Lombardo, Umberto. "Poles of Inaccessibility: A Calculation Algorithm for the Remotest Places on Earth." Scottish Geographical Journal 123 (2007): 227 – 233. Print.

[2] Garner, Will. "Shortest Distance from a Point to a Line." UCSD - Department of Mathematics . Web. 6 Dec. 2012. <http://math.ucsd.edu/~wgarner/math4c/derivations/distance/distptline.htm>.

[3] Wald, Abraham. "Statistical decision functions which minimize the maximum risk." Sequential tests of statistical hypotheses. New York: Annals of Mathematics, 1945. 265-280. Print.